\newif\ifHTML
\newif\ifAttachments
\HTMLfalse
\Attachmentstrue

\newcommand{\reforatt}[2]{%
  \attachfile{#1}~#2%
}%
\newcommand{\lb}{\\ &}
\newcommand{\wlb}{\\[1em] &}

\newcommand{\HorP}[2]{\ifHTML{}#1{}\else{}#2{}\fi}  

\newcommand{\htmlversion}{\url{http://markov-bases.de/models/K3N/K3N.html}}

\documentclass[a4paper,          
               12pt,             
               bibliography=totoc,         
               index=totoc,
               notitlepage]{scrartcl}
\usepackage[utf8]{inputenc}
\usepackage{amsmath}
\usepackage{amsthm}
\usepackage{amsfonts}
\usepackage{amssymb}

\usepackage{attachfile}
\usepackage{hyperref}
\usepackage{url}

\makeatletter
\@ifpackageloaded{tex4ht}{%
  \HTMLtrue%
  \Attachmentstrue%
  \renewcommand{\reforatt}[2]{%
    \href{#1}{#2}%
  }%
  \renewcommand{\lb}{}
  \renewcommand{\wlb}{\quad}
}{}
\makeatother
\usepackage{tikz}

\usepackage{graphicx}
\usepackage{color}

\usepackage{verbatim}

\newtheorem{lemma}{Lemma}

\newtheorem{thm}[lemma]{Theorem}

\theoremstyle{definition}

\theoremstyle{remark}

\newcommand{\ol}{\overline}

\newcommand{\Nb}{\mathbb{N}}

\newcommand{\Zb}{\mathbb{Z}}

\newcommand{\Fbf}{\mathbf{F}}

\newcommand{\Acal}{\mathcal{A}}
\newcommand{\Bcal}{\mathcal{B}}

\newcommand{\TFP}[1][\Acal]{\times_{#1}}

\allowdisplaybreaks[2]

\makeatletter

\newcounter{fourb@pos}

\newif\iffourb@negative

\def\@fourbdelim{,}
\def\@fourbsign{-}
\def\@fourb@plus{+}
\def\@fourb@minus{-}
\def\@fourb@one{1}

\newcommand{\fourb@printfield}{
  \iffourb@negative
    \ifx\fourb@field\@fourb@one
    -
    \else
    \overline{\fourb@field}
    \fi
  \else
    \ifx\fourb@field\@fourb@one
    +
    \else
    \fourb@field
    \fi
  \fi
}

\newcommand{\fourb@sep}{
  \fourb@printfield
  \gdef\fourb@field{}
  \fourb@negativefalse
  \stepcounter{fourb@pos}
  \def\fourb@separator{}
  \ifnum\value{fourb@pos}=3
  \def\fourb@separator{\\}
  \fi
  \ifnum\value{fourb@pos}=7
  \def\fourb@separator{\\}
  \fi
  \ifnum\value{fourb@pos}=9
  \def\fourb@separator{\endmatrix
    \\[1em]
    \matrix}
  \setcounter{fourb@pos}{1}
  \fi
  \ifnum\value{fourb@pos}=5
  \def\fourb@separator{\endmatrix
    &
    \matrix}
  \fi
  \fourb@separator
}

\newcommand{\fourb@step}[1]{
  \ifx\relax#1%
  \fourb@printfield
  \else%
     \if\@fourbdelim#1
        \fourb@sep%
     \else%
        \if\@fourbsign#1
           \fourb@negativetrue
        \else
           \xdef\fourb@field{\fourb@field#1}%
        \fi
     \fi
     \expandafter\fourb@step%
  \fi%
}

\newcommand{\fourb}[1]{%
  \setcounter{fourb@pos}{1}%
  \gdef\fourb@field{}
  \begin{matrix}
    \matrix
      \fourb@step #1\relax
    \endmatrix
  \end{matrix}
}

\makeatother

\author{
  \begin{minipage}{0.4\linewidth}
    \centering
    Johannes Rauh\\ \medskip\small
    MPI MIS\\
    Inselstraße 22\\
    04103 Leipzig\\
    jarauh@gmx.net
  \end{minipage}
  \begin{minipage}{0.4\linewidth}
    \centering
    Seth Sullivant\\ \medskip\small
    Department of Mathematics, NCSU\\
    Box 8205,\\
    Raleigh, NC 27695\\
    smsulli2@ncsu.edu 
  \end{minipage}}
\title{The Markov basis of $K_{3,N}$}

\begin{document}
\maketitle

\begin{abstract}
  This \HorP{website}{document} explains how to obtain a Markov basis of the graphical model of the complete bipartite
  graph~$K_{3,N}$ with binary nodes.  The computations illustrate the theory developed in~\cite{RS14:Lifting} that
  explains how to compute Markov bases of toric fiber products.

  \ifHTML\else%
  An HTML version of this document is available at \htmlversion. %
  \fi
\end{abstract}

\tableofcontents

\newcommand{\Kfh}[1][1]{\begin{tikzpicture}[scale=#1]  
    \path (0,2) coordinate (X1);
    \path (0.5,1) coordinate (X2);
    \path (1.5,1) coordinate (X3); 
    \path (0,0) coordinate (X4);
    \draw (X1) -- (X2) -- (X4) -- (X3) -- (X1);
    \draw (X2) -- (X3);
  \end{tikzpicture}}
\newcommand{\sKfh}{\Kfh[0.15]}  
\newcommand{\Kft}[1][1]{\begin{tikzpicture}[scale=#1]  
    \path (0,2) coordinate (X1);
    \path (0.5,1) coordinate (X2);
    \path (1.5,1) coordinate (X4);
    \path (0,0) coordinate (X3);
    \draw (X1) -- (X4) -- (X2);
    \draw (X3) -- (X4);
    \draw[fill] (X1) -- (X2) -- (X3) -- cycle;
  \end{tikzpicture}}
\newcommand{\sKft}{\Kft[0.15]}  
\newcommand{\Ts}[1][1]{\begin{tikzpicture}[scale=#1]  
    \path (0,2) coordinate (X1);
    \path (0,1) coordinate (X2);
    \path (1.5,1) coordinate (X4);
    \path (0,0) coordinate (X3);
    \draw (X1) -- (X4) -- (X2);
    \draw (X3) -- (X4);
  \end{tikzpicture}}
\newcommand{\sTs}{\Ts[0.15]}  

\section{Summary}
\label{sec:intro}

We compute Markov bases 
of the binary (i.e.~$d_{i}=2$) hierarchical model of the complete bipartite graph $K_{3,N}$.
\begin{thm}
  For any~$N$, the Markov degree of the binary hierarchical model of the complete bipartite graph $K_{3,N}$ is at most~$12$.

  The degree of the kernel Markov basis is at most~6, the degree of the PF Markov basis is~4 and the lifting defect is~2.
  Therefore, another bound on the degree of $K_{3,N}$ is~$4+2N$.
\end{thm}
The proof of this theorem will take up Sections~\ref{sec:Holes-K4t} to~\ref{sec:sKft-tableau} of this manuscript.  The
proof will also give an explicit description of a Markov basis of~$K_{3,N}$.  In Section~\ref{sec:comp-with-comp} we
compare our theoretical bounds with Markov bases that were obtained with \verb|4ti2|~\cite{4ti2}.  Both bounds are not
sharp for~$N\le 3$.

The proof relies heavily on the lifting machinery developped in~\cite{RS14:Lifting}.  All the notation and all the
notions are explained in detail in that manuscript.

The idea of the proof is to use the fact that $K_{3,N}$ is a toric fiber product of $N$ copies of the
three-star~$\sTs$.  The associated codimension zero product is a product of copies of the graph $\tilde K_{4}$ that
arises from $K_{4}$ by filling the triangle $\{1,2,3\}$ filled (\sKft).
The calculation is complicated by the fact that the marginal cone of $\sKft$ is not normal.
However, it turns out that all holes are vertices of the projected fibers.  This allows to treat the holes in a
systematic way by adding additional inequalities to the inequality description of the projected fibers.  We compute the
holes in Section~\ref{sec:Holes-K4t}.

The set of holes is described in Section~\ref{sec:Holes-K4t}.  This allows to describe the projected fibers and to
calculate a PF Markov basis (Section~\ref{sec:ineq-Kfh}).  The liftings are computed in
Section~\ref{sec:lifting-Kfh-tableau}, where it is also shown that the degree of the glued moves is at most~12.
Section~\ref{sec:sKft-tableau} presents a kernel Markov basis of degree six.  As it turns out, all moves from the kernel
Markov basis are redundant except for the quadratic moves.  The appendix contains the kernel Markov basis
(Appendix~\ref{sec:sKft-tensor}), the PF Markov basis (Appendix~\ref{sec:PF-MB-tensor}) and the lifts
(Appendix~\ref{sec:lifting-Kfh-tensor}) in tensor notation.



The results in this manuscript were obtained with the help of \verb|Normaliz|~\cite{Normaliz}, \verb|4ti2|~\cite{4ti2} and~\verb|Macaulay2|~\cite{M2}.
\ifHTML
This page contains links to files with code to reproduce the results.
\else
\ifAttachments
This PDF-file contains attached files with code to reproduce the results.  In PDF viewers that support attached files
they can usually be accessed by right clicks on the symbol \noattachfile{}.
\fi

An HTML-version of the results, which may be better suited for reading on a screen, is available at \htmlversion. %
\fi

\section{\texorpdfstring{The holes of $\tilde K_4$}{The holes}}
\label{sec:Holes-K4t}

In this section we study the set of holes of $\Nb\tilde K_4$.  The following lemma summarizes the most important properties:
\begin{lemma}
  \label{lem:holes-tK4}
  The set of holes of $\Nb\Bcal_{\tilde K_{4}}$ satisfies the following statements:
  \begin{enumerate}
  \item There are two fundamental holes~$h_{1},h_{2}$.
  \item There is a partition $\Bcal_{\tilde K_{4}}=\Bcal_{1}\dot\cup\Bcal_{2}$ of the rows of~$\Bcal_{\tilde K_{4}}$
    such that the set of holes is given by $(h_{1} + \Nb\Bcal_{1}) \dot\cup (h_{2} + \Nb\Bcal_{2})$.  The $\sTs$-margins of
    $\Bcal_{i}$ are linearly independent for $i=1,2$.
  \item There are linear functionals $l_{1},l_{2}:\Zb^{\Bcal_{\tilde K_{4}}}\to\Zb$ that satisfy the following:
    If a fiber $\Fbf(b)$ has a hole~$h\in(h_{i}+\Nb\Bcal_{i})$, then $l_{3-i}(v)>0=l_{3-i}(h)$.
  \end{enumerate}
\end{lemma}
The lemma follows from the observations made in the remainder of this section.
\ifAttachments
A Macaulay2-program that does the calculations in this section can be found
\reforatt{K3N-data/holesMonoid.m2}{here}.
\fi

A computation with \verb|Normaliz|
\ifAttachments
(\reforatt{K3N-data/K31codz.in}{input}/\reforatt{K3N-data/K31codz.out}{output})
\fi
shows that the Hilbert basis of the saturation of $\Nb\Bcal_{\tilde
  K_{4}}$ contains two additional vectors $h_{1},h_{2}$.  Both restrict to the all-ones hole $h$ on the
$K_{4}$-marginals (that is, on all pair marginals).  On the triangle-marginal, $h_{1}$ corresponds to XOR and $h_{2}$
corresponds to the opposite of XOR.  They are the only fundamental holes, in the nomenclature
of~\cite{HemmeckeTakemureYoshida09:Computing_holes_in_semigroups}, as the following computations show:
\begin{equation*}
  h_{1}+h_{2} = \Bcal_{\sKft}
  \begin{bmatrix}
    0000 \\
    0011 \\
    0101 \\
    0110 \\
    1000 \\
    1011 \\
    1101 \\
    1110
  \end{bmatrix},
  \qquad
  2h_{1} = \Bcal_{\sKft}
  \begin{bmatrix}
    0000 \\
    0001 \\
    0110 \\
    0111 \\
    1010 \\
    1011 \\
    1100 \\
    1101
  \end{bmatrix},
  \qquad
  2h_{2} = \Bcal_{\sKft}
  \begin{bmatrix}
    0010 \\
    0011 \\
    0100 \\
    0101 \\
    1000 \\
    1001 \\
    1110 \\
    1111
  \end{bmatrix}.
\end{equation*}

Consider one of the fundamental holes~$h_{i}$.  According
to~\cite{HemmeckeTakemureYoshida09:Computing_holes_in_semigroups}, we need to do the following:
\begin{enumerate}
\item Find the minimal non-negative solutions $(\lambda,\mu)$ of $h_{i} + \Bcal\lambda = \Bcal\mu$.
  This can be done, for example, using \verb|4ti2|'s command \verb|zsolve|.

  \ifAttachments
  Input:
  \reforatt{K3N-data/h1.mat}{\texttt{h1.mat}}/\reforatt{K3N-data/h1.mat}{\texttt{h2.mat}}, 
  \reforatt{K3N-data/h1.rhs}{\texttt{h1.rhs}}/\reforatt{K3N-data/h2.rhs}{\texttt{h2.rhs}},
  \ifHTML\else\\\phantom{Input: }\fi
  \reforatt{K3N-data/h1.sign}{\texttt{h1.sign}}/\reforatt{K3N-data/h2.sign}{\texttt{h2.sign}}.
  \\
  Output: \reforatt{K3N-data/h1.rhs}{\texttt{h1.zinhom}}/\reforatt{K3N-data/h2.rhs}{\texttt{h2.zinhom}}.
  \fi
\item Drop the $\mu$'s and interprete the $\lambda$'s (first half of the matrix) as the exponent vectors of monomials
  that generate a monomial ideal~$I$.
\item Compute the standard monomials.  In Macaulay2, this can be done using the command \verb|standardPairs|.
\end{enumerate}
The result is the following: Let
\begin{align*}
  \Bcal_{1}&=(b_{0000},b_{1100},b_{1010},b_{0110},b_{0001},b_{1101},b_{1011},b_{0111}),  \\
  \Bcal_{2}&=(b_{1000},b_{0100},b_{0010},b_{1110},b_{1001},b_{0101},b_{0011},b_{1111})
\end{align*}
($\Bcal_{1}$ corresponds to XOR on
the first three nodes, and $\Bcal_{2}$ corresponds to its opposite).  The holes derived from $h_{1}$ are of the form
$h_{1} + \Bcal_{1}\lambda$, where $\lambda\in\Nb^{8}$.
By symmetry, the holes derived from~$h_{2}$ are $h_{2} + \Bcal_{2}\mu$, where $\mu\in\Nb^{8}$.

Consider the two linear forms
\begin{align*}
  l_{1} &= y^{123}_{000} + y^{123}_{011} + y^{123}_{101} + y^{123}_{110}, \\
  l_{2} &= y^{123}_{001} + y^{123}_{010} + y^{123}_{100} + y^{123}_{111},
\end{align*}
where $y^{123}_{ijk}$ counts the coordinates where the $(123)$-marginal is equal to $ijk$ (that is, $l_{1}$ counts the
XOR-part of the triangle-margin, and $l_{2}$ counts the opposite XOR-part).  Then
\begin{equation*}
  l_{1}(h_{1} + \Bcal_{1}\lambda) > 0,
  \quad
  l_{2}(h_{1} + \Bcal_{1}\lambda) = 0,
  \quad
  l_{1}(h_{2} + \Bcal_{2}\mu) = 0,
  \quad
  l_{2}(h_{2} + \Bcal_{2}\mu) > 0.
\end{equation*}
Therefore, each hole is either derived from $h_{1}$ or from~$h_{2}$, 
that is, $(h_{1}+\Nb\Bcal_{1})\cap(h_{2}+\Nb\Bcal_{2})=\emptyset$.

Each set $\Bcal_{i}$ is linearly independent.  Hence different choices of the $\lambda$ (or $\mu$) give
different holes.  Moreover, also the vectors $\sTs$-marginals of $\Bcal_{i}$ are linearly independent for each~$i$.
Therefore, no two holes of the same type have the same $\sTs$-marginals.
Therefore, no projected fiber contains more than one hole of the same type. 
A projected fiber can have one hole of each type, though (for example, the holes $h_{1} + \Bcal_{1}(1,\dots,1)^{t}$ and
$h_{2} + \Bcal_{2}(1,\dots,1)^{t}$ have the same pair marginals and lie in the same fiber).

Let $h=h_{1}+\Bcal_{1}\lambda$ be a hole.  If $v$ belongs to the same projected fiber as~$h$, then
$v=\Bcal_{1}\lambda'+\Bcal_{2}\mu'$ with $\lambda',\mu'\in\Nb^{8}$.  As mentioned above, the $\sTs$-pair marginals of
$\Bcal_{1}$ are linearly independent.
Therefore, if $\mu'=0$, then, since
$h$ and $v$ have the same $\sTs$-marginals, $h=v$.  Hence, if $v\neq h$ is not a hole itself, $\mu'\neq 0$.  It follows
that $l_{2}(m)>0=l_{2}(h)$.

\section{The projected fiber Markov basis}
\label{sec:ineq-Kfh}

By Lemma~\ref{lem:holes-tK4}, every hole is a vertex of its projected fiber, supported either by $l_{1}$ or~$l_{2}$.
Thus we can do the following: We start with an inequality description of the semigroup $\Nb\tilde K_{4}$.  This gives us
a set of valid inequalities $D u \ge c$ for each projected fiber.  These inequalities are also valid for the holes.  We
augment the matrix $D$ by two additional rows corresponding to $l_{1}$ and~$l_{2}$ and denote the augmented matrix by~$D'$.  Then each projected fiber equals a solution set of linear inequalities of the form $D' u \ge c'$.
Therefore, any inequality Markov basis of~$D'$ can be used as a PF Markov basis.

Each projected fiber is a subset of~$\Zb^{8}$, with basis $e_{000},e_{001},\dots,e_{111}$.  Let $y_{000},\dots,y_{111}$
be the corresponding coordinates.  In a projected fiber, there are relations
\begin{align*}
  y_{011} &= y_{0}^{1} - y_{000} - y_{001} - y_{010}, \\
  y_{101} &= y_{0}^{2} - y_{000} - y_{001} - y_{100}, \\
  y_{110} &= y_{0}^{3} - y_{000} - y_{010} - y_{100}, \\
  y_{111} 
         & = 1 - y_{0}^{1} - y_{0}^{2} - y_{0}^{3} + 2 y_{000} + y_{001} + y_{010} + y_{100},
\end{align*}
where $y^{i}_{j}$ is the sum of those marginals where the $i$th entry equals $j$. %
Since the projected fiber is four-dimensional, these are all relations.  An independent set of coordinates is given by
$y_{000},y_{001},y_{010},y_{100}$ (that is, those coordinates with at most one one).  According to \verb|Normaliz|%
\ifAttachments
\ (\reforatt{K3N-data/K31codz.in}{input}/\reforatt{K3N-data/K31codz.out}{output})%
\fi%
, they satisfy inequalities of the form
\begin{gather*}
  ? \le y_{000}, \quad ? \le y_{001}, \quad ? \le y_{010}, \quad ? \le y_{100}, \\
  ? \le y_{000} + y_{001} \le ?, \quad
  ? \le y_{000} + y_{010} \le ?, \quad
  ? \le y_{000} + y_{100} \le ?, \\
  y_{000} + y_{001} + y_{010} \le ?, \quad
  y_{000} + y_{001} + y_{100} \le ?, \quad
  y_{000} + y_{010} + y_{100} \le ?, \\
  ? \le 2 y_{000} + y_{001} + y_{010} + y_{100}.
\end{gather*}
To get rid of the holes, we need to add the inequalities with linear parts $l_{1},l_{2}$.  In coordinates, they take the form
\begin{equation*}
  ? \le 2 (y_{000} + y_{001} + y_{010} + y_{100}) \le ?. 
\end{equation*}

The columns of the corresponding matrix $D$ spans a lattice, and \verb|4ti2| computes its Markov basis%
\ifAttachments
\ (\reforatt{K3N-data/K31codz-iMBp.lat}{input}/\reforatt{K3N-data/K31codz-iMBp.mar}{output})%
\fi.
Since the first four rows of $D$ form a unit matrix, $D$ is easy to invert: The first four coordinates of the Markov
basis of $D$ give the PF Markov basis in the coordinates $y_{000},y_{001},y_{010},y_{100}$.
In tableau notation, PF Markov basis consists of the 16 moves, that are (up to symmetry; that is, up to a permutation of
the columns) of the form
\begin{align*}
  \begin{bmatrix}
    0 & 0 & a \\
    1 & 1 & b
  \end{bmatrix}
  &-
  \begin{bmatrix}
    0 & 1 & a \\
    1 & 0 & b
  \end{bmatrix},
  &
  \begin{bmatrix}
    0 & 0 & 0 \\
    0 & 0 & 1 \\
    1 & 1 & 0 \\
    1 & 1 & 1
  \end{bmatrix}
  &-
  \begin{bmatrix}
    0 & 1 & 0 \\
    0 & 1 & 1 \\
    1 & 0 & 0 \\
    1 & 0 & 1
  \end{bmatrix},
  &
  \begin{bmatrix}
    0 & 0 & 0 \\
    0 & 1 & 1 \\
    1 & 0 & 1 \\
    1 & 1 & 0
  \end{bmatrix}
  &-
  \begin{bmatrix}
    0 & 0 & 1 \\
    0 & 1 & 0 \\
    1 & 0 & 0 \\
    1 & 1 & 1
  \end{bmatrix}.
\end{align*}

\section{Lifting the PF Markov basis}
\label{sec:lifting-Kfh-tableau}

In this section we compute the lifts.  We use the algorithm described in~\cite{RS14:Lifting} 
to compute the lifting as an inequality Markov basis.
Analyzing the results we find that the maximal degree of a glue of lifts is bounded by~12. %
\ifAttachments %
The file \reforatt{K3N-data/lift.m2}{lift.m2} contains Macaulay2 code that does the calculations in this section.  It
makes use of further routines from \reforatt{K3N-data/M2routines.m2}{M2routines.m2}. %
\fi

First consider
\begin{equation*}
  g =
  \begin{bmatrix}
    0 & 0 & a \\
    1 & 1 & b
  \end{bmatrix}
  -
  \begin{bmatrix}
    0 & 1 & a \\
    1 & 0 & b
  \end{bmatrix}.
\end{equation*}
If $b=a$, then $g$ lifts to
\begin{align*}
  \begin{bmatrix}
    0 & 0 & a & c \\
    1 & 1 & a & c
  \end{bmatrix}
  -
  \begin{bmatrix}
    0 & 1 & a & c \\
    1 & 0 & a & c
  \end{bmatrix},
  \qquad &
  \begin{bmatrix}
    0 & 0 & a & c \\
    1 & 1 & a & \ol{c} \\
    1 & d & \ol{a} & c \\
    0 & d & \ol{a} & \ol{c}
  \end{bmatrix}
  -
  \begin{bmatrix}
    1 & 0 & a & c \\
    0 & 1 & a & \ol{c} \\
    0 & d & \ol{a} & c \\
    1 & d & \ol{a} & \ol{c}
  \end{bmatrix}
  \lb
  \text{ and }\quad
  \begin{bmatrix}
    0 & 0 & a & c \\
    1 & 1 & a & \ol{c} \\
    d & 1 & \ol{a} & c \\
    d & 0 & \ol{a} & \ol{c}
  \end{bmatrix}
  -
  \begin{bmatrix}
    0 & 1 & a & c \\
    1 & 0 & a & \ol{c} \\
    d & 0 & \ol{a} & c \\
    d & 1 & \ol{a} & \ol{c}
  \end{bmatrix}.
\end{align*}
The lifting defect is at most two.  Any move $\tilde m$ that arises by gluing lifts of $g$ satisfies
\begin{equation*}
  \xi(\tilde m^{+}) - g^{+} \le
  \begin{bmatrix}
    1 & 0 & \ol a \\
    0 & 0 & \ol a \\
    1 & 1 & \ol a \\
    0 & 1 & \ol a
  \end{bmatrix}.
\end{equation*}
Therefore, $\deg(\tilde m)\le 6$.

If $b=\ol{a}$, then $g$ lifts to
\begin{align*}
  \begin{bmatrix}
    0 & 0 & a & c \\
    1 & 1 & \ol{a} & c
  \end{bmatrix}
  -
  \begin{bmatrix}
    0 & 1 & a & c \\
    1 & 0 & \ol{a} & c
  \end{bmatrix},
  \quad &
  \begin{bmatrix}
    0 & 0 & a & c \\
    1 & 1 & \ol{a} & \ol{c} \\
    0 & 0 & \ol{a} & \ol{c} \\
    0 & 1 & \ol{a} & c
  \end{bmatrix}
  -
  \begin{bmatrix}
    0 & 1 & a & c \\
    1 & 0 & \ol{a} & \ol{c} \\
    0 & 1 & \ol{a} & \ol{c} \\
    0 & 0 & \ol{a} & c
  \end{bmatrix},
  \lb
  \text{and }\quad
  \begin{bmatrix}
    0 & 0 & a & c \\
    1 & 1 & \ol{a} & \ol{c} \\
    1 & 0 & a & \ol{c} \\
    1 & 1 & a & c
  \end{bmatrix}
  -
  \begin{bmatrix}
    0 & 1 & a & c \\
    1 & 0 & \ol{a} & \ol{c} \\
    1 & 1 & a & \ol{c} \\
    1 & 0 & a & c
  \end{bmatrix}.
\end{align*}
The lifting defect is at most two.  Any move $\tilde m$ that arises by gluing lifts of $g$ satisfies
\begin{equation*}
  \xi(\tilde m^{+}) - g^{+} \le
  \begin{bmatrix}
    0 & 0 & \ol a \\
    0 & 1 & \ol a \\
    1 & 0 & a \\
    1 & 1 & a
  \end{bmatrix}.
\end{equation*}
Therefore, $\deg(\tilde m)\le 6$.

For
\begin{equation*}
  g =
  \begin{bmatrix}
    0 & 0 & 0 \\
    0 & 0 & 1 \\
    1 & 1 & 0 \\
    1 & 1 & 1
  \end{bmatrix}
  -
  \begin{bmatrix}
    0 & 1 & 0 \\
    0 & 1 & 1 \\
    1 & 0 & 0 \\
    1 & 0 & 1
  \end{bmatrix}
\end{equation*}
the lifts are (up to symmetry, exchanging the first two columns and state switching)
\begin{align*}
  \begin{bmatrix}
    0 & 0 & 0 & a \\
    0 & 0 & 1 & b \\
    1 & 1 & 0 & a \\
    1 & 1 & 1 & b
  \end{bmatrix}
  - 
  \begin{bmatrix}
    0 & 1 & 0 & a \\
    0 & 1 & 1 & b \\
    1 & 0 & 0 & a \\
    1 & 0 & 1 & b
  \end{bmatrix},
  \quad 
  \begin{bmatrix}
    0 & 0 & 0 & a \\
    0 & 0 & 1 & b \\
    1 & 1 & 0 & b \\
    1 & 1 & 1 & a
  \end{bmatrix}
  - &
  \begin{bmatrix}
    0 & 1 & 0 & a \\
    0 & 1 & 1 & b \\
    1 & 0 & 0 & b \\
    1 & 0 & 1 & a
  \end{bmatrix},
  \wlb
  \begin{bmatrix}
    0 & 0 & 0 & a \\
    0 & 0 & 1 & \ol{a} \\
    1 & 1 & 0 & a \\
    1 & 1 & 1 & a \\
    0 & 0 & 0 & a \\
    1 & 0 & 0 & \ol{a}
  \end{bmatrix}
  - 
  \begin{bmatrix}
    1 & 0 & 0 & a \\
    1 & 0 & 1 & \ol{a} \\
    0 & 1 & 0 & a \\
    0 & 1 & 1 & a \\
    1 & 0 & 0 & a \\
    0 & 0 & 0 & \ol{a}
  \end{bmatrix}.
\end{align*}
The lifting defect is at most two.  Any move $\tilde m$ that arises by gluing lifts of $g$ satisfies
\begin{equation*}
  \xi(\tilde m^{+}) - g^{+} \le
  \begin{bmatrix}
    0 & 0 & 0 \\
    1 & 0 & 0 \\
    0 & 1 & 0 \\
    1 & 1 & 0 \\
    0 & 0 & 1 \\
    1 & 0 & 1 \\
    0 & 1 & 1 \\
    1 & 1 & 1
  \end{bmatrix}.
\end{equation*}
Therefore, $\deg(\tilde m)\le 12$.

For
\begin{equation*}
  \begin{bmatrix}
    0 & 0 & 0 \\
    0 & 1 & 1 \\
    1 & 0 & 1 \\
    1 & 1 & 0
  \end{bmatrix}
  -
  \begin{bmatrix}
    1 & 0 & 0 \\
    1 & 1 & 1 \\
    0 & 0 & 1 \\
    0 & 1 & 0
  \end{bmatrix}  
\end{equation*}
the lifts are (up to symmetry, exchanging the first two columns and state switching)
\begin{align*}
  \begin{bmatrix}
    0 & 0 & 0 & a \\
    0 & 1 & 1 & b \\
    1 & 0 & 1 & b \\
    1 & 1 & 0 & a
  \end{bmatrix}
  &-
  \begin{bmatrix}
    1 & 0 & 0 & a \\
    1 & 1 & 1 & b \\
    0 & 0 & 1 & b \\
    0 & 1 & 0 & a
  \end{bmatrix},
  &
  \begin{bmatrix}
    0 & 0 & 0 & a \\
    0 & 1 & 1 & b \\
    1 & 0 & 1 & a \\
    1 & 1 & 0 & b
  \end{bmatrix}
  &-
  \begin{bmatrix}
    1 & 0 & 0 & a \\
    1 & 1 & 1 & b \\
    0 & 0 & 1 & a \\
    0 & 1 & 0 & b
  \end{bmatrix},
  \\
  \begin{bmatrix}
    0 & 0 & 0 & a \\
    0 & 1 & 1 & \ol{a} \\
    1 & 0 & 1 & a \\
    1 & 1 & 0 & a \\
    0 & 0 & 0 & a \\
    1 & 0 & 0 & \ol{a}
  \end{bmatrix}
  &-
  \begin{bmatrix}
    1 & 0 & 0 & a \\
    1 & 1 & 1 & \ol{a} \\
    0 & 0 & 1 & b \\
    0 & 1 & 0 & a \\
    1 & 0 & 0 & a \\
    0 & 0 & 0 & \ol{a}
  \end{bmatrix},
  &
  \begin{bmatrix}
    0 & 0 & 0 & \ol{a} \\
    0 & 1 & 1 & a \\
    1 & 0 & 1 & a \\
    1 & 1 & 0 & a \\
    0 & 0 & 0 & \ol{a} \\
    1 & 0 & 0 & \ol{a}
  \end{bmatrix}
  &-
  \begin{bmatrix}
    1 & 0 & 0 & \ol{a} \\
    1 & 1 & 1 & a \\
    0 & 0 & 1 & a \\
    0 & 1 & 0 & a \\
    1 & 0 & 0 & \ol{a} \\
    0 & 0 & 0 & \ol{a}
  \end{bmatrix}
\end{align*}
The lifting defect is at most two.  Any move $\tilde m$ that arises by gluing lifts of $g$ satisfies
\begin{equation*}
  \xi(\tilde m^{+}) - g^{+} \le
  \begin{bmatrix}
    0 & 0 & 0 \\
    1 & 0 & 0 \\
    0 & 1 & 0 \\
    1 & 1 & 0 \\
    0 & 0 & 1 \\
    1 & 0 & 1 \\
    0 & 1 & 1 \\
    1 & 1 & 1
  \end{bmatrix}.
\end{equation*}
Therefore, $\deg(\tilde m)\le 12$.

\section{The kernel Markov basis in tableau notation} 
\label{sec:sKft-tableau}

The Markov basis of $\sKft$ computed by \verb|4ti2| has 20 elements of degrees four and six%
\ifAttachments
\ (\reforatt{K3N-data/K31codz.mat}{input}/\reforatt{K3N-data/K31codz.mar}{output})%
\fi.  In tableau notation, it
consists of the following moves (up to symmetry, involving permuting the first three columns):
\begin{equation*}
  \begin{bmatrix} 
    a000 \\ b011 \\ a101 \\ b110
  \end{bmatrix}
  -
  \begin{bmatrix}
    a001 \\ b010 \\ a100 \\ b111
  \end{bmatrix}
  , 
  \begin{bmatrix} 
    abc0 \\ abc0 \\ \ol a\ol b\ol c0 \\ \ol abc1 \\ a\ol bc1 \\ ab\ol c1
  \end{bmatrix}
  -
  \begin{bmatrix}
    abc1 \\ abc1 \\ \ol a\ol b\ol c1 \\ \ol abc0 \\ a\ol bc0 \\ ab\ol c0
  \end{bmatrix}.
\end{equation*}

Recall that the kernel Markov basis consists of quadratic moves of the form
\begin{equation*}
  \begin{bmatrix} 
    abcD^{\phantom{\prime}}E^{\phantom{\prime}} \\ abcD'E'
  \end{bmatrix}
  -
  \begin{bmatrix}
    abcD^{\phantom{\prime}}E' \\ abcD'E^{\phantom{\prime}}
  \end{bmatrix}  
\end{equation*}
and lifts that can be constructed from the moves in the Markov basis of $\sKft$ 
by adding constant columns in the tableau notation.  The lifted moves have degree 4 and~6, and so the kernel Markov
basis has degree~6.

Let us show that all lifted moves from the kernel Markov basis are redundant; that is: We only need the quadratic moves
from the kernel Markov basis.  In fact, consider, for example, a lift $m$ of degree six.  The glues of lifts of the PF
Markov basis contain the quadratic moves of the form
\begin{equation*}
  \begin{bmatrix}
    a^{\phantom{\prime}}b^{\phantom{\prime}}c^{\phantom{\prime}}d^{\phantom{\prime}}E \\ a'b'c'd'E
  \end{bmatrix}
  -
  \begin{bmatrix}
    a'b^{\phantom{\prime}}c^{\phantom{\prime}}d^{\phantom{\prime}}E \\ a^{\phantom{\prime}}b'c'd'E
  \end{bmatrix},
\end{equation*}
and so on.  Applying such quadratic moves reduces $m$ to zero:
\begin{equation*}
  \begin{bmatrix}
    abc0E \\ abc0E \\ \ol a\ol b\ol c0E \\ \ol abc1E \\ a\ol bc1E \\ ab\ol c1E
  \end{bmatrix}
  \to
  \begin{bmatrix}
    \ol abc0E \\ abc0E \\ a\ol b\ol c0E \\ \ol abc1E \\ a\ol bc1E \\ ab\ol c1E
  \end{bmatrix}
  \to
  \begin{bmatrix}
    \ol abc0E \\ a\ol bc0E \\ ab\ol c0E \\ \ol abc1E \\ a\ol bc1E \\ ab\ol c1E
  \end{bmatrix}
  \to
  \begin{bmatrix}
    \ol abc0E \\ a\ol bc0E \\ ab\ol c0E \\ abc1E \\ a\ol bc1E \\ \ol ab\ol c1E
  \end{bmatrix}
  \to
  \begin{bmatrix}
    \ol abc0E \\ a\ol bc0E \\ ab\ol c0E \\ abc1E \\ abc1E \\ \ol a\ol b\ol c1E
  \end{bmatrix}
  =
  \begin{bmatrix}
    abc1E \\ abc1E \\ \ol a\ol b\ol c1E \\ \ol abc0E \\ a\ol bc0E \\ ab\ol c0E
  \end{bmatrix}.
\end{equation*}
The quartic moves of the kernel Markov basis can be reduced similarly.

\section{Comparison with computational results}
\label{sec:comp-with-comp}

For $N\le 3$, the Markov basis of $K_{3,N}$ can be computed (within reasonable time) using \verb|4ti2|.  The Markov
degrees are:
\begin{center}
  \begin{tabular}{l|ccc}
    $N$    & 1 & 2 & 3 \\
    $\deg$ & 2 & 4 & 6
  \end{tabular}
\end{center}
These three computed degrees are much smaller than the theoretical bound of~$\min\{4+2 N,12\}$.

$K_{3,1}$ is a tree; hence the Markov degree is two.  The additional moves obtained by lifting the PF Markov basis are
not necessary.  For a zero-fold toric fiber product it is no wonder that our bound is far from being tight.


A Markov basis of $K_{3,2}$ was computed in~\cite{KRS11:Positive_Margins_and_Binomial_Walks} by interpreting $K_{3,2}$
as a TFP of three two-stars.  Here, we interprete these moves from the viewpoint of our above computations.  The Markov
basis of~\cite{KRS11:Positive_Margins_and_Binomial_Walks} consists of two kinds of quadrics and quartics.  First, there
are the codimension-zero quadrics.  Second, the quadrics of the two three-stars glue together and yield further quadrics.
Observe that the two interpretations of $K_{3,2}$ interchange the roles of the codimension-zero quadrics and the glued
quadrics: The codimension-zero quadrics of $K_{1,2}\TFP[\Acal'] K_{1,2} \TFP[\Acal'] K_{1,2}$ correspond to the glued quadrics of
$K_{3,1}\TFP K_{3,1}$, and vice versa.

The quartics are of the form
\begin{equation*}
  Q =
  \begin{bmatrix}
    0 & a_{00} & b_{00} & 0 & 0 \\
    1 & a_{01} & b_{01} & 0 & 1 \\
    1 & a_{10} & b_{10} & 1 & 0 \\
    0 & a_{11} & b_{11} & 1 & 1
  \end{bmatrix}
  -
  \begin{bmatrix}
    1 & a_{00} & b_{00} & 0 & 0 \\
    0 & a_{01} & b_{01} & 0 & 1 \\
    0 & a_{10} & b_{10} & 1 & 0 \\
    1 & a_{11} & b_{11} & 1 & 1
  \end{bmatrix}
\end{equation*}
for some $l_{ij},m_{ij}\in\{0,1\}$, modulo a permutation of the first three rows.  
All of these quadrics are glues of, say, $m$ and~$m'$.  If either $(a_{00},b_{00})=(a_{01},b_{01})$ or
$(a_{10},b_{10})=(a_{11},b_{11})$, then $m$ is a quadric; that is
\begin{equation*}
  m =
  \begin{bmatrix}
    1 & a_{10} & b_{10} & 1 \\
    0 & a_{11} & b_{11} & 1
  \end{bmatrix}
  -
  \begin{bmatrix}
    0 & a_{10} & b_{10} & 1 \\
    1 & a_{11} & b_{11} & 1
  \end{bmatrix},
  \text{ or }
  m =
  \begin{bmatrix}
    1 & a_{00} & b_{00} & 1 \\
    0 & a_{01} & b_{01} & 1
  \end{bmatrix}
  -
  \begin{bmatrix}
    0 & a_{00} & b_{00} & 1 \\
    1 & a_{01} & b_{01} & 1
  \end{bmatrix}.
\end{equation*}
Otherwise, $m$ is the quartic
\begin{equation*}
  m =
  \begin{bmatrix}
    0 & a_{00} & b_{00} & 0 \\
    1 & a_{01} & b_{01} & 0 \\
    1 & a_{10} & b_{10} & 1 \\
    0 & a_{11} & b_{11} & 1
  \end{bmatrix}
  -
  \begin{bmatrix}
    1 & a_{00} & b_{00} & 0 \\
    0 & a_{01} & b_{01} & 0 \\
    0 & a_{10} & b_{10} & 1 \\
    1 & a_{11} & b_{11} & 1
  \end{bmatrix}.
\end{equation*}
That is, $m$ is a lift of
\begin{equation*}
  g =
  \begin{bmatrix}
    0 & a_{00} & b_{00} \\
    1 & a_{01} & b_{01} \\
    1 & a_{10} & b_{10} \\
    0 & a_{11} & b_{11}
  \end{bmatrix}
  -
  \begin{bmatrix}
    1 & a_{00} & b_{00} \\
    0 & a_{01} & b_{01} \\
    0 & a_{10} & b_{10} \\
    1 & a_{11} & b_{11}
  \end{bmatrix}.  
\end{equation*}
A $g$ of this form may be quadratic or quartic, depending on the $(a_{ij},b_{ij})$.  One can see that projecting the
Markov basis of $K_{3,2}$ gives the full PF Markov basis.
Therefore, the reason that the bound is not sharp is not that the PF Markov basis is too large, but that not all glues
are necessary.

\ifHTML\else





\fi

\appendix

\section{Intermediate results in tensor notation}
\label{sec:tensor-notation}

In this appendix we present some of the above results in tensor notation.  These results are straightforward
translations of the output of~\verb|4ti2|, without factoring out any symmetry or other structures, and so they should be
considered as intermediate steps on the way from the raw \verb|4ti2|-output towards the form summarized above.

The tensor notation is as follows: Any state $x\in\{0,1\}^{4}$ corresponds to a binary string~$x_{4}x_{3}x_{2}x_{1}$,
which can be considered as the binary representation of a natural number.  The three least-significant bits
$x_{3}x_{2}x_{1}$ correspond to the first group of nodes in~$K_{3,N}$.
The states appear at the following positions:
\begin{equation*}
  \fourb{0000\ , 0001, 0010\ , 0011, 0100\ , 0101, 0110\ , 0111,
         1000\ , 1001, 1010\ , 1011, 1100\ , 1101, 1110\ , 1111}.
\end{equation*}
Entries of modulus one are replaced by their signs, and negative numbers are indicated by a $\overline{\text{bar}}$.


\subsection{The kernel Markov basis in tensor notation} 
\label{sec:sKft-tensor}

The Markov basis of $\sKft$ computed by \verb|4ti2| has 20 elements of degrees four and six:
\begin{align*}
  \fourb{1,  0, -1,  0, -1,  0,  1,  0, -1,  0,  1,  0,  1,  0, -1,  0}, &&  
  \fourb{0,  1,  0, -1,  0, -1,  0,  1,  0, -1,  0,  1,  0,  1,  0, -1}, &&  
  \fourb{1, -1,  0,  0, -1,  1,  0,  0, -1,  1,  0,  0,  1, -1,  0,  0}, &&  
  \fourb{0,  0,  1, -1,  0,  0, -1,  1,  0,  0, -1,  1,  0,  0,  1, -1}, &&  
  \fourb{1, -1, -1,  1,  0,  0,  0,  0, -1,  1,  1, -1,  0,  0,  0,  0}, &&  
  \fourb{0,  0,  0,  0,  1, -1, -1,  1,  0,  0,  0,  0, -1,  1,  1, -1}, \\[1em]  
  \fourb{1,  0, -1,  0,  0, -1,  0,  1, -1,  0,  1,  0,  0,  1,  0, -1}, &&  
  \fourb{0,  1,  0, -1, -1,  0,  1,  0,  0, -1,  0,  1,  1,  0, -1,  0}, &&  
  \fourb{1, -1,  0,  0,  0,  0, -1,  1, -1,  1,  0,  0,  0,  0,  1, -1}, &&  
  \fourb{0,  0,  1, -1, -1,  1,  0,  0,  0,  0, -1,  1,  1, -1,  0,  0}, &&  
  \fourb{1,  0,  0, -1, -1,  0,  0,  1, -1,  0,  0,  1,  1,  0,  0, -1}, &&  
  \fourb{0,  1, -1,  0,  0, -1,  1,  0,  0, -1,  1,  0,  0,  1, -1,  0}, \\[1em]  
  &&
  \fourb{2, -1, -1,  0, -1,  0,  0,  1, -2,  1,  1,  0,  1,  0,  0, -1}, &&  
  \fourb{1, -2,  0,  1,  0,  1, -1,  0, -1,  2,  0, -1,  0, -1,  1,  0}, &&  
  \fourb{1,  0, -2,  1,  0, -1,  1,  0, -1,  0,  2, -1,  0,  1, -1,  0}, &&  
  \fourb{0,  1,  1, -2, -1,  0,  0,  1,  0, -1, -1,  2,  1,  0,  0, -1}, \\[1em]  
  &&
  \fourb{1,  0,  0, -1, -2,  1,  1,  0, -1,  0,  0,  1,  2, -1, -1,  0}, &&  
  \fourb{0,  1, -1,  0,  1, -2,  0,  1,  0, -1,  1,  0, -1,  2,  0, -1}, &&  
  \fourb{0,  1, -1,  0, -1,  0,  2, -1,  0, -1,  1,  0,  1,  0, -2,  1}, &&  
  \fourb{1,  0,  0, -1,  0, -1, -1,  2, -1,  0,  0,  1,  0,  1,  1, -2} . 
\end{align*}

\subsection{The projected fiber Markov basis in tensor notation}
\label{sec:PF-MB-tensor}

The PF Markov basis consists of the 16 $2\times2\times2$-tableaus
\begin{align*}
  &\begin{matrix}
    +- & 00 \\
    -+ & 00
  \end{matrix},
  &&\begin{matrix}
    00 & +- \\
    00 & -+
  \end{matrix},
  &&\begin{matrix}
    +0 & -0 \\
    -0 & +0
  \end{matrix},
  &&\begin{matrix}
    0+ & 0- \\
    0- & 0+
  \end{matrix},
  &&\begin{matrix}
    +- & -+ \\
    00 & 00
  \end{matrix},
  &&\begin{matrix}
    00 & 00 \\
    +- & -+
  \end{matrix},
  \\[1em]
  &\begin{matrix}
    +0 & 0- \\
    -0 & 0+
  \end{matrix},
  &&\begin{matrix}
    0+ & -0 \\
    0- & +0
  \end{matrix},
  &&\begin{matrix}
    +- & 00 \\
    00 & -+
  \end{matrix},
  &&\begin{matrix}
    00 & -+ \\
    +- & 00
  \end{matrix},
  &&\begin{matrix}
    +0 & -0 \\
    0- & 0+
  \end{matrix},
  &&\begin{matrix}
    0- & 0+ \\
    +0 & -0
  \end{matrix},
\end{align*}
\begin{align*}
  &\begin{matrix}
    +- & +- \\
    -+ & -+
  \end{matrix},
  &&\begin{matrix}
    +- & -+ \\
    +- & -+
  \end{matrix},
  &&\begin{matrix}
    ++ & -- \\
    -- & ++
  \end{matrix},
  \\[1em]
  &&&\begin{matrix}
    +- & -+ \\
    -+ & +-
  \end{matrix}.
\end{align*}

\subsection{Lifting in tensor notation}
\label{sec:lifting-Kfh-tensor}


For
\begin{equation*}
  \fourb{1,-1,-1,1,0,0,0,0}
\end{equation*}
there are ten lifts:
\begin{align*}
  &&
  \fourb{1, -1, -1, 1, 0, 0, 0, 0, 0, 0, 0, 0, 0, 0, 0, 0},   &&  
  \fourb{0, 0, 0, 0, 0, 0, 0, 0, 1, -1, -1, 1, 0, 0, 0, 0},   \\[1em]  
  \fourb{1, -1, 0, 0, 0, 0, -1, 1, 0, 0, -1, 1, 0, 0, 1, -1}, &&  
  \fourb{1, -1, 0, 0, -1, 1, 0, 0, 0, 0, -1, 1, 1, -1, 0, 0}, &&  
  \fourb{0, 0, 1, -1, -1, 1, 0, 0, -1, 1, 0, 0, 1, -1, 0, 0}, &&  
  \fourb{0, 0, 1, -1, 0, 0, -1, 1, -1, 1, 0, 0, 0, 0, 1, -1}, \\[1em]  
  \fourb{0, 1, 0, -1, -1, 0, 1, 0, -1, 0, 1, 0, 1, 0, -1, 0}, &&  
  \fourb{1, 0, -1, 0, -1, 0, 1, 0, 0, -1, 0, 1, 1, 0, -1, 0}, &&  
  \fourb{0, 1, 0, -1, 0, -1, 0, 1, -1, 0, 1, 0, 0, 1, 0, -1}, &&  
  \fourb{1, 0, -1, 0, 0, -1, 0, 1, 0, -1, 0, 1, 0, 1, 0, -1}. 
\end{align*}

For
\begin{equation*}
  \fourb{1,0,-1,0,0,-1,0,1}
\end{equation*}
there are six lifts:
\begin{align*}
  &&
  \fourb{1, 0, -1, 0, 0, -1, 0, 1, 0, 0, 0, 0, 0, 0, 0, 0},   &&
  \fourb{0, 0, 0, 0, 0, 0, 0, 0, 1, 0, -1, 0, 0, -1, 0, 1},   \\[1em]
  \fourb{0, 0, 0, 0, 1, -1, -1, 1, 1, 0, -1, 0, -1, 0, 1, 0}, &&
  \fourb{0, 1, 0, -1, 0, -1, 0, 1, 1, -1, -1, 1, 0, 0, 0, 0}, &&
  \fourb{1, -1, -1, 1, 0, 0, 0, 0, 0, 1, 0, -1, 0, -1, 0, 1}, &&
  \fourb{1, 0, -1, 0, -1, 0, 1, 0, 0, 0, 0, 0, 1, -1, -1, 1}.
\end{align*}

For
\begin{equation*}
  \fourb{1,-1,-1,1,1,-1,-1,1}
\end{equation*}
there are 21 lifts:
\begin{align*}
  \fourb{0, 1, 0, -1, -1, 0, 1, 0, -1, 0, 1, 0, 0, 1, 0, -1} , 
\end{align*}
\begin{align*}
  \fourb{1, -1, -1, 1, 1, -1, -1, 1, 0, 0, 0, 0, 0, 0, 0, 0} ,  &&  
  \fourb{0, 0, 0, 0, 0, 0, 0, 0, 1, -1, -1, 1, 1, -1, -1, 1} ,  && 
  \fourb{1, -1, -1, 1, 0, 0, 0, 0, 0, 0, 0, 0, 1, -1, -1, 1} ,  &&  
  \fourb{0, 0, 0, 0, 1, -1, -1, 1, 1, -1, -1, 1, 0, 0, 0, 0} ,  \\[1em]  
  \fourb{2, -2, -1, 1, 0, 0, -1, 1, -1, 1, 0, 0, 1, -1, 0, 0},  &&  
  \fourb{1, -1, -2, 2, 1, -1, 0, 0, 0, 0, 1, -1, 0, 0, -1, 1},  &&  
  \fourb{2, -1, -2, 1, 0, -1, 0, 1, -1, 0, 1, 0, 1, 0, -1, 0},  &&  
  \fourb{1, -2, -1, 2, 1, 0, -1, 0, 0, 1, 0, -1, 0, -1, 0, 1},  \\[1em]  
  \fourb{0, 0, 1, -1, -2, 2, 1, -1, -1, 1, 0, 0, 1, -1, 0, 0},  &&  
  \fourb{1, -1, 0, 0, 1, -1, -2, 2, 0, 0, -1, 1, 0, 0, 1, -1},  &&  
  \fourb{0, 1, 0, -1, -2, 1, 2, -1, -1, 0, 1, 0, 1, 0, -1, 0},  &&  
  \fourb{1, 0, -1, 0, 1, -2, -1, 2, 0, -1, 0, 1, 0, 1, 0, -1},  \\[1em]  
  \fourb{1, -1, 0, 0, -1, 1, 0, 0, -2, 2, 1, -1, 0, 0, 1, -1},  &&  
  \fourb{0, 0, 1, -1, 0, 0, -1, 1, 1, -1, -2, 2, 1, -1, 0, 0},  &&  
  \fourb{1, 0, -1, 0, -1, 0, 1, 0, -2, 1, 2, -1, 0, 1, 0, -1},  &&  
  \fourb{0, 1, 0, -1, 0, -1, 0, 1, 1, -2, -1, 2, 1, 0, -1, 0},  \\[1em]  
  \fourb{1, -1, 0, 0, -1, 1, 0, 0, 0, 0, -1, 1, 2, -2, -1, 1},  &&  
  \fourb{0, 0, 1, -1, 0, 0, -1, 1, -1, 1, 0, 0, -1, 1, 2, -2},  &&  
  \fourb{1, 0, -1, 0, -1, 0, 1, 0, 0, -1, 0, 1, 2, -1, -2, 1},  &&  
  \fourb{0, 1, 0, -1, 0, -1, 0, 1, -1, 0, 1, 0, -1, 2, 1, -2}. 
\end{align*}

For
\begin{equation*}
  \fourb{1,-1,-1,1,-1,1,1,-1}
\end{equation*}
there are 40 lifts:
\begin{align*}
  \fourb{1, -1, -1, 1, -1, 1, 1, -1, 0, 0, 0, 0, 0, 0, 0, 0},   &&  
  \fourb{0, 0, 0, 0, 0, 0, 0, 0, 1, -1, -1, 1, -1, 1, 1, -1},   &&  
  \fourb{1, -1, -1, 1, 0, 0, 0, 0, 0, 0, 0, 0, -1, 1, 1, -1},   &&  
  \fourb{0, 0, 0, 0, 1, -1, -1, 1, -1, 1, 1, -1, 0, 0, 0, 0},   \\[1em]  
  \fourb{1, -1, 0, 0, -1, 1, 0, 0, 0, 0, -1, 1, 0, 0, 1, -1},   &&  
  \fourb{0, 0, 1, -1, 0, 0, -1, 1, -1, 1, 0, 0, 1, -1, 0, 0},   &&  
  \fourb{1, 0, -1, 0, -1, 0, 1, 0, 0, -1, 0, 1, 0, 1, 0, -1},   &&  
  \fourb{0, 1, 0, -1, 0, -1, 0, 1, -1, 0, 1, 0, 1, 0, -1, 0},   \\[1em]  
  \fourb{2, -2, -1, 1, -1, 1, 0, 0, -1, 1, 0, 0, 0, 0, 1, -1},   &&  
  \fourb{1, -1, -2, 2, 0, 0, 1, -1, 0, 0, 1, -1, -1, 1, 0, 0},   &&  
  \fourb{2, -1, -2, 1, -1, 0, 1, 0, -1, 0, 1, 0, 0, 1, 0, -1},   &&  
  \fourb{1, -2, -1, 2, 0, 1, 0, -1, 0, 1, 0, -1, -1, 0, 1, 0},   \\[1em]  
  \fourb{1, -1, 0, 0, -2, 2, 1, -1, 0, 0, -1, 1, 1, -1, 0, 0},   &&  
  \fourb{0, 0, 1, -1, 1, -1, -2, 2, -1, 1, 0, 0, 0, 0, 1, -1},   &&  
  \fourb{1, 0, -1, 0, -2, 1, 2, -1, 0, -1, 0, 1, 1, 0, -1, 0},   &&  
  \fourb{0, 1, 0, -1, 1, -2, -1, 2, -1, 0, 1, 0, 0, 1, 0, -1},   \\[1em]  
  \fourb{1, -1, 0, 0, 0, 0, -1, 1, -2, 2, 1, -1, 1, -1, 0, 0},   &&  
  \fourb{0, 0, 1, -1, -1, 1, 0, 0, 1, -1, -2, 2, 0, 0, 1, -1},   &&  
  \fourb{1, 0, -1, 0, 0, -1, 0, 1, -2, 1, 2, -1, 1, 0, -1, 0},   &&  
  \fourb{0, 1, 0, -1, -1, 0, 1, 0, 1, -2, -1, 2, 0, 1, 0, -1},   \\[1em]  
  \fourb{0, 0, 1, -1, -1, 1, 0, 0, -1, 1, 0, 0, 2, -2, -1, 1},   &&  
  \fourb{1, -1, 0, 0, 0, 0, -1, 1, 0, 0, -1, 1, -1, 1, 2, -2},   &&  
  \fourb{0, 1, 0, -1, -1, 0, 1, 0, -1, 0, 1, 0, 2, -1, -2, 1},   &&  
  \fourb{1, 0, -1, 0, 0, -1, 0, 1, 0, -1, 0, 1, -1, 2, 1, -2},   \\[1em]  
  \fourb{2, -1, -1, 0, -2, 1, 1, 0, -1, 0, 0, 1, 1, 0, 0, -1},   &&  
  \fourb{1, -2, 0, 1, -1, 2, 0, -1, 0, 1, -1, 0, 0, -1, 1, 0},   &&  
  \fourb{1, 0, -2, 1, -1, 0, 2, -1, 0, -1, 1, 0, 0, 1, -1, 0},   &&  
  \fourb{0, 1, 1, -2, 0, -1, -1, 2, -1, 0, 0, 1, 1, 0, 0, -1},   \\[1em]  
  \fourb{1, 0, 0, -1, -1, 0, 0, 1, -2, 1, 1, 0, 2, -1, -1, 0},   &&  
  \fourb{0, 1, -1, 0, 0, -1, 1, 0, 1, -2, 0, 1, -1, 2, 0, -1},   &&  
  \fourb{0, 1, -1, 0, 0, -1, 1, 0, -1, 0, 2, -1, 1, 0, -2, 1},   &&  
  \fourb{1, 0, 0, -1, -1, 0, 0, 1, 0, -1, -1, 2, 0, 1, 1, -2},   \\[1em]  
  \fourb{2, -1, -1, 0, -1, 0, 0, 1, -1, 0, 0, 1, 0, 1, 1, -2},   &&  
  \fourb{1, -2, 0, 1, 0, 1, -1, 0, 0, 1, -1, 0, -1, 0, 2, -1},   &&  
  \fourb{0, 1, 1, -2, -1, 0, 0, 1, -1, 0, 0, 1, 2, -1, -1, 0},   &&  
  \fourb{1, 0, -2, 1, 0, -1, 1, 0, 0, -1, 1, 0, -1, 2, 0, -1},   \\[1em]  
  \fourb{1, 0, 0, -1, -2, 1, 1, 0, 0, -1, -1, 2, 1, 0, 0, -1},   &&  
  \fourb{0, 1, -1, 0, 1, -2, 0, 1, -1, 0, 2, -1, 0, 1, -1, 0},   &&  
  \fourb{0, 1, -1, 0, -1, 0, 2, -1, 1, -2, 0, 1, 0, 1, -1, 0},   &&  
  \fourb{1, 0, 0, -1, 0, -1, -1, 2, -2, 1, 1, 0, 1, 0, 0, -1},   \\[1em]  
\end{align*}


\bibliographystyle{IEEEtranS}
\bibliography{HighCodim}

\end{document}